\documentclass[a4paper,11pt]{article}

\usepackage{latexsym}
\usepackage{amssymb}        
\newtheorem{thm}{Theorem}[section]

\newtheorem{cor}[thm]{Corollary}
\newtheorem{prop}[thm]{Proposition}

\newtheorem{rem}[thm]{Remark}

\input epsf

\title{On generating series of complementary planar trees}
\author{Roland Bacher
}
\date{}
\begin{document}
\maketitle

{\it Abstract\footnote{Math. class: 05A15, 05C05, 06A10. Keywords: Integer sequence, generating function, inversion of power series, planar tree, spin model}:} 
We generalize and reprove an identity
of Parker and Loday. It states that certain pairs of
generating series associated
to pairs of labelled rooted planar trees are mutually inverse
under composition.

\section{Introduction} In \cite{CSV} Carlitz, Scoville and Vaughan
consider finite words (in a finite alphabet $\mathcal A$) such that
all pairs of consecutive letters belong to a fixed subset 
$L\subset {\mathcal A}\times {\mathcal A}$. They show (Theorems 6.8 and 
7.3 of \cite{CSV}) that suitably defined
pairs of signed generating series counting such words associated 
to $L\subset {\mathcal A}\times {\mathcal A}$ and to its complementary set 
${\overline L}={\mathcal A}\times 
{\mathcal A}\setminus L$ are each others inverse.
Their result was generalized in the first part of Parker's thesis \cite{P} 
who showed an analogous 
result for suitable classes of finite trees having labelled vertices.
Loday in \cite {L}, motivated by questions concerning 
combinatorial realisations of operads, rediscovered Parkers result
toghether with a different proof, based on homological arguments.
 
This paper presents a further generalisation of Parkers and Lodays
result.

A typical example of our identity can be described as follows:
Associate to two complex matrices
$$M_1=\left(\begin{array}{cc}a&b\\c&d\end{array}\right),\ 
M_2=\left(\begin{array}{cc}\alpha&\beta\\\gamma&\delta\end{array}\right)
$$
the following two systems of algebraic equations
$$\left\lbrace\begin{array}{lcl}
g_1&=&(-X+ag_1+bg_2)(-X+\alpha g_1+\beta g_2)\\
g_2&=&(-X+cg_1+dg_2)(-X+\gamma g_1+\delta g_2)\\
g&=&-X+g_1+g_2\end{array}\right.$$
and
$$\left\lbrace\begin{array}{lcl}
\tilde g_1&=&(X-(1-a)\tilde g_1-(1-b)\tilde g_2)(X-(1-\alpha)\tilde g_1-(1-
\beta)\tilde g_2)\\
\tilde g_2&=&(X-(1-c)\tilde g_1-(1-d)\tilde g_2)(X-(1-\gamma)\tilde  g_1-
(1-\delta)\tilde  g_2)\\
\tilde g&=&-X+\tilde g_1+\tilde g_2\end{array}\right.\ .$$
Choosing continuous determinations satisfying
$g=g_1=g_2=\tilde g=\tilde g_1=\tilde g_2=0$ at
$X=0$ we get holomorphic
functions $g=g(X),\ \tilde g=\tilde g(X)$ for $X$ 
in an open neighbourhood of $0\in{\mathbb C}$. We have now 
$$g(\tilde g(X))=X$$ for all $X$ in a small open 
disc centered at $0$. 

This result holds of course formally for the corresponding
generating series and can 
be verified by computing for instance a minimal polynomial
$P(u,v)=\sum_{i,j} p_{i,j}u^iv^j$ for $g$ (i.e. satisfying in particular
$P(g(X),X)=0$) and checking that we have $P(X,\tilde g(X))=0$. 
Since this identity is algebraic, the field ${\mathbb C}$ can be
replaced by an arbitrary commutative ring.

The sequel of this paper is 
organized as follows: The next section states the main result
in purely algebraic terms over a not necessarily 
commutative ring. Parkers and Lodays result corresponds to
the special case of matrices with coefficients
in $\{0,1\}$.
Our main result removes the restriction on the coefficients. 
It is also somewhat easier to state (at least in a commutative setting)
since it avoids combinatorial descriptions.
It follows from our formulation that all involved generating 
functions are algebraic in a commutative setting and over a finite 
alphabet.
Section 3 fixes notations concerning trees.
Section 4 describes spin models on trees and recasts our main result
using partition functions of spin models. Section 5 proves
the main result using a spin model on grafted trees (called \lq\lq graftings''
in \cite{L}). The proof avoids homological arguments and is thus
in some sense more elementary (although
perhaps not simpler) than the proof of \cite{L}. 
Section 6 describes briefly a further generalisation involving 
arbitrary (not necessarily regular) finite trees which appears already
in Parkers work.
Section 7 is a digression 
generalizing the notion of grafted trees. Section 8 contains
the computations for example (i) of \cite{L}. We display the
defining polynomial of the relevant (algebraic) generating function and 
discuss briefly its asymptotics.

\section{Main result}
Consider a (not necessarily finite) alphabet ${\mathcal A}$ 
and a (not necessarily commutative)
associative ring $R$ having a unit $1$. We denote by $Y_{\mathcal A}$
a set of (non-commutative) variables indexed by elements $\alpha\in {\mathcal
A}$ and by $X$ a supplementary (non-commutative) variable.
We denote by $\tilde R[[X]]={\mathbf Z}[[R,Y_{\mathcal A},R]]\otimes_{
\mathbf Z} {\mathbf Z}$ the ring of formal power series. An element
of $\tilde R[[X]]$ is a (generally infinite) sum of monomials of the form
$$\pm r_1Z_1r_2Z_2\cdots r_lZ_lr_{l+1}$$
with $r_i\in R$ and $Z_i\in \{X\}\cup Y_{\mathcal A}$. 
Let $k\geq 2$ be a natural integer and let $M_1,\dots,M_k$
be a set of matrices with rows and columns indexed by ${\mathcal A}$
and coefficients
$M_j(\alpha,\beta)\in R$ for $\alpha,\beta\in {\mathcal A}$. 
For $\alpha\in{\mathcal A}$, let $g_{\alpha}\in \tilde R[[X]]$ 
be the power series
$g_\alpha=Y_\alpha X^k+\hbox{ (terms of higher order in }X)$ which satisfies
$$g_\alpha=Y_\alpha(X-(M_1\ V)_\alpha)(X-(M_2\ V)_\alpha)\cdots
(X-(M_k\ V)_\alpha)$$
where $V$ is the column vector with coordinates 
$g_\beta,\ \beta\in{\mathcal A}$ and
where $(M_j\ V)_\alpha=\sum_{\beta\in {\mathcal A}} M_j(\alpha,\beta)g_\beta$ 
denotes the 
$\alpha-$th coordinate of the matrix-product $M_j\ V$.
Set $g=-X+\sum_{\alpha\in{\mathcal A}} g_\alpha$. 

Remark that our notation is slightly misleading: $g,g_1,\dots,g_l$ are
power series in $X$ defining \lq\lq functions'' of $X$. 
The letter $X$ stands of course for the power series at $X=0$ of
the identity function $X\longmapsto X$. Remark also that the requirements 
$k\geq 2$ and $g_\alpha=Y_\alpha X^k+\hbox{ terms of higher order in }X$ ensure
that $g_\alpha$ is well-defined: Its $n-$th coefficient involving $X$
$n$ times depends only on
coefficients with indices $\leq n-k+1$ of $g_\beta,\ \beta\in{_mathcal A}$ 
involving at most $n-k+1$ occurences of $X$.  

Define \lq\lq complementary'' matrices 
$\tilde M_1,\dots,\tilde M_k$ with coefficients
$\tilde M_j(\alpha,\beta)=1-M_j(\alpha,\beta)$ by
setting $\tilde M_j=J-M_j$ where $J$ is the all $1$ matrix
with rows and columns indexed by $\mathcal A$. 
For $\alpha\in{\mathcal A}$ we introduce \lq\lq complementary'' functions 
$\tilde g_\alpha=Y_\alpha(-X)^k+
\hbox{ terms of higher order in }X\in \tilde R[[X]]$ satisfying
$$\tilde g_\alpha=Y_\alpha(-X+(\tilde M_1\ \tilde V)_\alpha)\cdots
(-X+(\tilde M_k\ \tilde V)_\alpha)$$
where $\tilde V$ is the column vector with coordinates
$\tilde g_\beta,\ \beta\in{\mathcal A}$ and where
$(\tilde M_j\ \tilde V)_\alpha=\sum_{\beta\in{\mathcal A}} 
\tilde M_j(\alpha,\beta)\tilde g_\beta$. 
We define now $\tilde g=-X+\sum_{j=1}^l \tilde g_j$. 

Given two formal power series $f,h\in\tilde R[[X]]$ such that 
every monomial of $h$ is at least of degree $1$ in $X$, the composition
$f\circ_X h$ of $f$ with $h$ is defined as the formal power series
obtained by 
replacing every occurence of $X^k,\ k=1,2,3,\dots$ in every monomial of
$f$ by the series $h^k$.

The main result of this paper can now be stated as follows:

\begin{thm} \label{main} We have (formally) 
$$g\circ_X\tilde g=\tilde g\circ_X g=X\ .$$
\end{thm}

\begin{rem} (i) In a commutative setting with $R={\mathbf C}$ the field 
of complex numbers, a finite alphabet $\mathcal A$, and $Y_\alpha\in
{\mathbf C}$ for $\alpha\in {\mathcal A}$,
the functions 
$g_\alpha,\tilde g_\alpha$ and thus also $g$ and $\tilde g$ 
are holomorphic determinations of algebraic functions in an open 
neighbourhood of $0$.

\ \ (ii) The above definitions of $g_\alpha$ and $\tilde g_\alpha$
are well-suited for iterative computations of 
the power-series $g_\alpha,\tilde g_\alpha$ by 
\lq\lq bootstrapping'' (see e.g. section 5.4 of \cite{O}). 
Indeed, given $g_{\alpha,n},\ \alpha\in{\mathcal A}$ such that 
$g_\alpha-g_{\alpha,n}=O(X^{1+(n+1)(k-1)})$ we have 
$$g_\alpha-Y_\alpha(t-(M_1\ V_n)_\alpha)\cdots(t-(M_k\ V_n)_\alpha)=
O(X^{1+(n+2)(k-1)})$$
where $V_n$ is the column vector with coordinates 
$g_{\alpha,n},\ \alpha\in{\mathcal A}$. An analogous result holds
of course for $\tilde g_1,\dots,\tilde g_l$.

\ \ (iii) Changing the grading and considering also all variables 
$Y_\alpha,\ \alpha\in{\mathcal A}$ as beeing of degree $1$, 
Theorem \ref{main} holds also for $k=1$ (as a special case, one obtains
the main result of \cite{CSV}). In this case, $g\circ_X\tilde g$
boils down to a simple product in $\tilde R$ since $g$ and $\tilde g$
are of the form $c X$ with $c\in R[[(Y_\alpha)_{\alpha\in{\mathcal A}}]]$.

\ \ (iv) The choice of an integer $k\geq 2$ (or $k=1$)
corresponds to the case of $k-$regular trees. This restriction
can be removed: Section 6 describes a further generalization (already 
contained in \cite{P}), corresponding to arbitrary finite trees.
\end{rem}

\section{Trees}

A {\it tree} is a connected graph without cycles. A {\it rooted} tree
contains a marked vertex $r$, called the {\it root}. In particular,
a rooted tree is non-empty. A rooted tree $T=\{r\}$ reduced to 
its root is {\it trivial}. The edges of a rooted tree are canonically 
oriented by requiring the root to be the unique source of the directed
tree obtained by orienting all edges away from 
the root. We write $e=(\alpha(e),\omega(e))$ for the edge
$e$ oriented from $\alpha(e)$ to $\omega(e)$.
Given an edge $e=(\alpha(e),\omega(e))$
we call $\omega(e)$ a {\it son} of $\alpha(e)$ and $\alpha(e)$ the 
{\it father} of $\omega(e)$. Each vertex $v\not=r\in T$ 
except the root has a unique 
father. The set of vertices sharing a common father is a {\it brotherhood}.
A {\it leaf} is a vertex without sons. A vertex having at least
one son is {\it interior}. The {\it level} of a vertex 
is its (combinatorial) distance to the root. A rooted tree is {\it planar}
if every brotherhood is totally ordered. Every brotherhood of a 
{\it locally finite} rooted tree is finite . A rooted tree is thus locally
finite if and only if the set of vertices of given level $N$ is finite for 
all $N\in{\mathbf N}$. 
A rooted tree is {\it $k-$regular} if every non-empty
brotherhood of strictly positive level contains exactly $k$ vertices.
A rooted tree is finite if and only if 
the set of vertices of level $N$ is finite for all $N$ and
empty for $N$ large enough. Given a vertex $v\in T$, we denote by $T(v)$
the subtree of $T$ rooted at $v$ which is defined by considering 
all vertices $w\in T$ such that $v$ belongs to the unique geodesic path
joining the root $r\in T$ to $w$ (the vertices of $T(v)$ are thus 
$v$ and its descendants). We call such a subtree {\it maximal}.
A {\it principal} subtree of a non-trivial tree $T$ is a maximal
subtree $T(s_i)$ rooted at a son $s_i$ of the root $r \in T$.
A non-trivial $k-$regular tree has thus exactly $k$ principal 
subtrees.

Unless otherwise stated, a ($k-$regular) tree will henceforth always denote
a ($k-$regular) finite rooted planar tree.

Every tree $T$ can be represented by a plane tree of the halfplane 
${\mathbf R}\times{\mathbf R}_{\geq 0}$ such that a totally ordered
brotherhood $s_1<s_2<\dots$ of level $n\geq 0$ corresponds to vertices 
$(x_1,n),(x_2,n),\dots,\ x_1<x_2<\dots$.

A vertex $v\in T$ of level $n$ has a unique recursively defined {\it address} 
$$a(v)=a_0\ a_1\ \dots\ a_{n-1}\in \{1,2,\dots\}^n$$
where $a_{n-1}$ is the number of elements $\leq v$
in the brotherhood of $v$ and where $a_0\ \dots\ a_{n-2}$ is the
address of the (unique) father of $v$. The address of the 
root $r\in T$ is empty.
The lexicographic order on addresses orders the set of vertices of a 
planar rooted tree completely.

\vskip0.5cm

\centerline{\epsfysize3.5cm\epsfbox{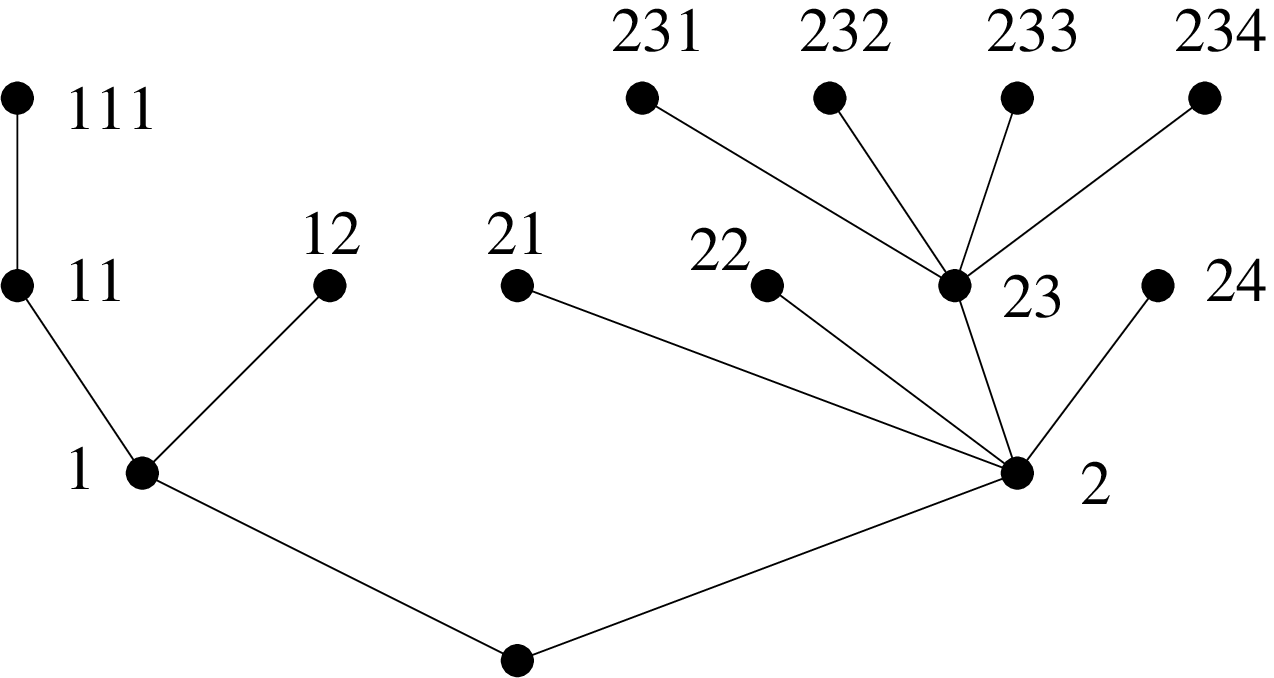}}
\centerline{Figure 1: A (finite rooted planar) tree with addressed
vertices.} 
\medskip

Figure 1 illustrates the notions of this section. It displays
a tree with 14 vertices
$$\emptyset<1<11<111<12<2<21<22<23<231<232<233<234<24\ ,$$
completely ordered by their addresses. The root corresponds 
to the lowest vertex
having an empty address. The oriented edge $e$ with vertices
$23$ and $232$ starts at $\alpha(e)=(23)$ and ends at the leave 
$\omega(e)=232$. The vertex $23$ is thus the father of $232$.
The brotherhood of $232$ are the four vertices $231,\ 232,\ 233$ 
and $234$ at level $3$ having $\alpha(e)=23$ as their common father. 
This tree is of course not regular 
since the interior vertices 1 and 2 (for instance) have
respectively $2$ and $4$ sons.
The tree of Figure 1 has two principal subtrees. Two vertices 
with a non-empty adress belong to a common principal subtree if
and only if the first letter of their adresses coincides.

\section{Non-commutative spin models on rooted planar trees}

Let $T=(E,V)$ be a (rooted planar) tree with edges $E$
and vertices $V$. We denote by $V^\circ$ the set of interior vertices 
(having at least one son) of $T$ and by $E'\subset E$ the subset of 
leafless edges (i.e. $e\in E'$ if $\omega(e)\in V^\circ$).
We define a {\it spin-model} $(T,w)$ on $T$ by considering  
a (not necessarily) finite set or alphabet $\mathcal A$ of 
{\it spins} and a {\it weight-function} 
$$w:E'\times {\mathcal A}\times {\mathcal A}\longrightarrow R$$
with values in a (not necessarily commutative) ring $R$ containing $1$. 
Setting $w_e(\alpha,\beta)=w(e,\alpha,\beta)$,
a weight function can be identified with an application
$$E'\longrightarrow M_{\mathcal A}(R)$$
where $M_{\mathcal A}(R)$ denotes the set of matrices 
with coefficients in $R$
and rows and columns indexed by $\mathcal A$. We call
$w_e\in M_{\mathcal A}(R)$ the {\it weight-matrix} of $e\in E'$.
We introduce furthermore 
(non-commutative) variables $Y_\alpha$ for $\alpha\in {\mathcal A}$
and a supplementary (non-commutative) variable $X$.

Given a spin-model $(T,w)$, its {\it complementary spin-model} is defined as
$(T,\tilde w)$ where $\tilde w_e(\alpha,\beta)=
1-w_e(\alpha,\beta)$ for $e\in E'$ and $\alpha,\beta\in {\mathcal A}$.

A colouring
$$\varphi:V^\circ\longrightarrow {\mathcal A}$$
of all interior vertices in $T$ by elements of $\mathcal A$ is a {\it state}.
Its {\it energy} $f_T(\varphi)\in R[X,(Y_\alpha)_{\alpha\in{\mathcal A}}]$ 
is recursively defined as follows:
If $T=\{r\}$ is trivial then $f_T(\varphi)=X$. Otherwise, we
set
$$f_T(\varphi)=Y_{\varphi(r)}\left(w(1,\varphi)f_{T_1}(\varphi_1)\right)
\left(w(2,\varphi)f_{T_2}(\varphi_2)\right)\cdots$$
where $T_1,T_2,\dots$ are the
principal subtrees associated to the linearly ordered sons 
$s_1<s_2<\dots$ of the root $r\in T$, where 
$$w(i,\varphi)=\left\{\begin{array}{ll}
w((r,s_i),\varphi(r),\varphi(s_i))&\hbox{if }s_i\in V^\circ\\
1&\hbox{otherwise}\end{array}\right.$$
and where $f_{T_i}(\varphi)$ is the energy of
the spinmodel on the principal 
subtree $T_i$ defined by the $i-$th son $s_i$ of the root $r\in T$
with weights and state obtained by restriction.

In a commutative setting, this boils down to 
$$f(\varphi)=X^{\sharp(V\setminus V^\circ)}\prod_{v\in V^\circ}Y_{\varphi(v)}
\prod_{e\in E'} w_e(\varphi(\alpha(e)),\varphi(\omega(e)))$$
which is the familiar definition used in statistical physics.

The total sum 
$$Z=Z(T)=\sum_{\varphi\in {\mathcal A}^{V^\circ}}f(\varphi)$$
of energies over all states is the {\it partition function}. We denote by 
$$Z_\alpha=Z_\alpha(T)=\sum_{\varphi\in {\mathcal A}^{V^\circ},\ 
\varphi(r)=\alpha}
f(\varphi)$$
the {\it restricted} partition function
obtained by computing the total energy of all states with prescribed colour 
$\varphi(r)=\alpha$ on the root. We have obviously
$$Z=\sum_{\alpha\in {\mathcal A}}Z_\alpha\ $$
if $T$ is non-trivial. For $T=\{r\}$ trivial, we have $Z=X$ and $Z_\alpha=0$
for $\alpha\in {\mathcal A}$.

The partition function of a tree $T$ can be computed as 
follows: We have $Z(T)=X$ for $T=\{r\}$ the trivial tree. 
Otherwise, consider the principal subtrees $T_1,T_2,\dots$
associated to the $k$ linearly ordered sons $s_1<s_2<\dots<s_k$ of $r$.
For $\alpha\in{\mathcal A}$ denote by $Z_\alpha(T_j)$ 
the obvious restricted partition functions  
of the principal subtree $T_j$ with colour $\alpha$ on its root $s_j$. 
We denote by $e_j
=\{r=\alpha(e_j),s_j=\omega(e_j)\}$ the oriented
edge joining the root $r$ of $T$ to its $j-$th son $s_j=\omega(e_j)$.

The definition of the restricted partition function implies then
easily the following result:

\begin{prop} If $T\not=\{r\}$ is non-trivial, we have
$$Z_\alpha(T)=Y_\alpha C_1(\alpha)C_2(\alpha)\cdots C_k(\alpha)$$
where $T_1,\dots,T_k$ are the principal subtrees associated to the
$k$ linearly ordered sons $s_1<s_2<\dots<s_k$ of the root $r\in T$ and where
$$C_k(\alpha)=\left\{\begin{array}{ll}\sum_{\beta\in{\mathcal A}}
w((r,s_i),\alpha,\beta)Z_\beta(T_i)&\hbox{if }T_i\not=\{s_i\}\\
X&\hbox{otherwise .}\end{array}\right.$$
\end{prop}

An edge $e\in E$ of a $k-$regular tree is of type $j$ if the 
extremity $\omega(e)$ (having level $m+1$)
of $e$ has an address
$$a_0\ a_1\ \dots\ a_m=a_0\ a_1\ \dots a_{m-1}\ j$$
ending with $1\leq j=a_m\leq k$. The endpoint $\omega(e)$ of a type 
$j$ edge is thus the $j-$th element in its (totally ordered)
brotherhood.

Given a (finite) $k-$regular tree $T$ and $k$ weight-matrices 
$M_1,\dots,M_k$ (with indices in ${\mathcal A}\times {\mathcal A}$) 
we consider the spin model with spins
$\mathcal A$ and weight function $w_e(\alpha,\beta)=M_j(\alpha,\beta)$
on leafless edges of type $j$.

\vskip0.5cm

\centerline{\epsfysize3.5cm\epsfbox{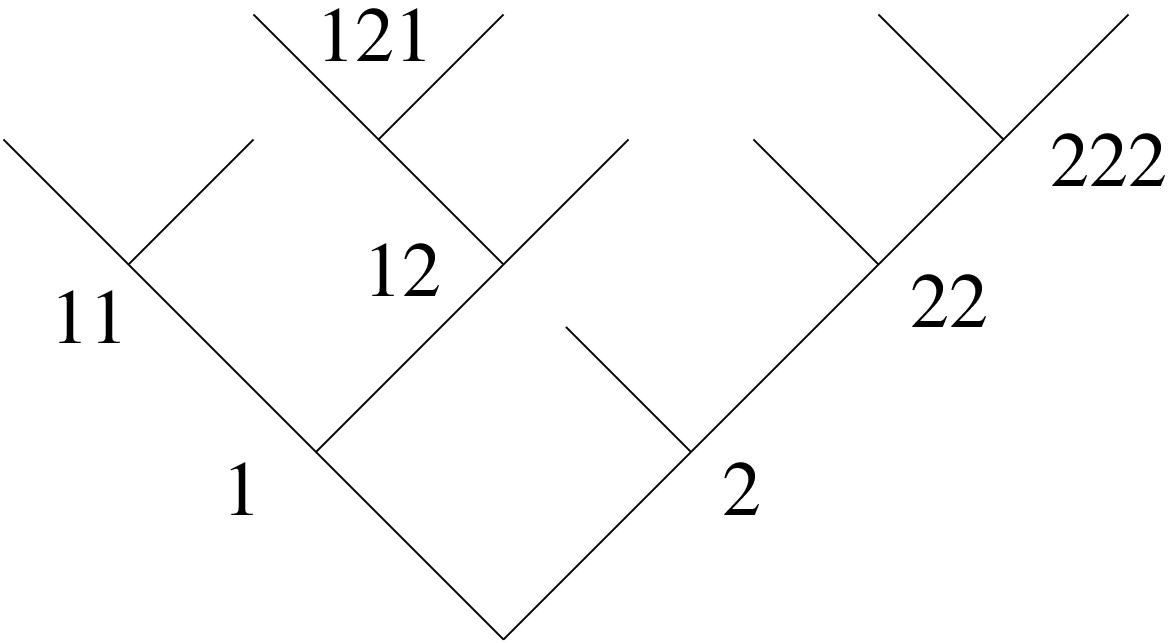}}
\centerline{Figure 2: A $2-$regular tree.} 
\medskip

{\bf Example.} We consider the commutative spin model on the $2-$regular
tree $T$ of Figure 2 with ${\mathcal A}=\{1,2\}$ and
$$M_1=\left(\begin{array}{rr}1&1\\1&2\end{array}\right),\ 
M_2=(M_1)^{-1}=\left(\begin{array}{rr}2&-1\\-1&1\end{array}\right)\ .$$
For the sake of simplicity, we set $X=Y_1=Y_2=1$.
Given an interior vertex $v$ of $T$, we denote by
$Z_*(v)$ the vector $\left(\begin{array}{c}Z_1(T(v))\\Z_2(T(v))
\end{array}\right)$ where $T(v)$ denotes the maximal
subtree of $T$ with root $v$ (obtained by considering the connected component
of $v$ in $T\setminus a(v)$ where $a(v)$ is the father of $v\not=
r$). Since all sons of $11,\ 121,\ 222$ are leaves, we have
$$Z_*(11)=Z_*(121)=Z_*(222)=\left(\begin{array}{c}1\\1
\end{array}\right)$$
where we denote a vertex by its address. We get now 
$$Z_*(12)=\left(\begin{array}{rr}1&1\\1&2\end{array}\right)
Z_*(121)=
\left(\begin{array}{rr}1&1\\1&2\end{array}\right)
\left(\begin{array}{c}1\\1\end{array}\right)=
\left(\begin{array}{c}2\\3\end{array}\right)\ ,$$
$$Z_*(22)=\left(\begin{array}{rr}2&-1\\-1&1\end{array}\right)
Z_*(222)=
\left(\begin{array}{rr}2&-1\\-1&1\end{array}\right)
\left(\begin{array}{c}1\\1\end{array}\right)=
\left(\begin{array}{c}1\\0\end{array}\right)$$
and
$$Z_*(2)=\left(\begin{array}{rr}2&-1\\-1&1\end{array}\right)
Z_*(22)=
\left(\begin{array}{rr}2&-1\\-1&1\end{array}\right)
\left(\begin{array}{c}1\\0\end{array}\right)=
\left(\begin{array}{c}2\\-1\end{array}\right)\ .$$

Denoting by $U*W$ the Hadamard product (\lq\lq
student's vector product'') $\left(\begin{array}{c}\alpha\alpha'\\
\beta\beta'\end{array}\right)$ of two vectors $\left(\begin{array}{c}\alpha\\
\beta\end{array}\right),\left(\begin{array}{c}\alpha'\\
\beta'\end{array}\right)$ we have then
$$\begin{array}{rcl}
Z_*(1)&=&\left(\left(\begin{array}{rr}1&1\\1&2\end{array}\right)
Z_*(11)\right)*
\left(\left(\begin{array}{rr}2&-1\\-1&1\end{array}\right)
Z_*(12)\right)\\
&=&
\left(\left(\begin{array}{rr}1&1\\1&2\end{array}\right)
\left(\begin{array}{c}1\\1\end{array}\right)\right)*
\left(\left(\begin{array}{rr}2&-1\\-1&1\end{array}\right)
\left(\begin{array}{c}2\\3\end{array}\right)\right)\\
&=&
\left(\begin{array}{c}2\\3\end{array}\right)*
\left(\begin{array}{c}1\\1\end{array}\right)=
\left(\begin{array}{c}2\\3\end{array}\right)
\end{array}$$
and 
$$\begin{array}{rcl}
Z_*(\emptyset)&=&\left(\left(\begin{array}{rr}1&1\\1&2\end{array}\right)
Z_*(1)\right)*
\left(\left(\begin{array}{rr}2&-1\\-1&1\end{array}\right)
Z_*(2)\right)\\
&=&
\left(\left(\begin{array}{rr}1&1\\1&2\end{array}\right)
\left(\begin{array}{c}2\\3\end{array}\right)\right)*
\left(\left(\begin{array}{rr}2&-1\\-1&1\end{array}\right)
\left(\begin{array}{c}2\\-1\end{array}\right)\right)\\
&=&
\left(\begin{array}{c}5\\8\end{array}\right)*
\left(\begin{array}{c}5\\-3\end{array}\right)=
\left(\begin{array}{c}25\\-24\end{array}\right)\ .
\end{array}$$
This yields the partition function
$$Z(T)=Z_1(T)+Z_2(T)=25-24=1$$
for the tree of Figure 2 with respect to the weights
$M_1=\left(\begin{array}{cc}1&1\\1&2\end{array}\right)$ and $M_2=M_1^{-1}$
associated to leafless edges indicating $NW$ (type 1), respectively
$NE$ (type 2).

Given a spinset (alphabet) $\mathcal A$ and pairs of complementary
matrices $M_1,\tilde M_1=J-M_1,\dots,M_k,\tilde M_k=J-
M_k$ (where $J$ denotes the all $1$ matrix with row and column-indices in
${\mathcal A}$), we consider the associated complementary spin models on 
$k-$regular trees with partition function $Z(T),\ \tilde Z(T)$ as above.

Denoting by ${\mathcal T}_k$ the set of all (finite, rooted, planar) 
$k-$regular trees, 
we introduce the signed formal generating series
$$Z({\mathcal T}_k)=-\sum_{T\in {\mathcal T}_k}(-1)^{d^\circ(T)}\ 
Z(T)$$
and
$$\tilde Z({\mathcal T}_k)=\sum_{T\in {\mathcal T}_k}(-1)^{d(T)}
\tilde Z(T)$$
where $d^\circ(T)=\sharp(V^\circ)$ denotes the number of interior vertices and
$d(T)=\sharp(V\setminus V^\circ)$ 
the number of leaves of a $k-$regular tree $T\in{\mathcal T}_k$.

Let us moreover consider the restricted generating series
$$Z_\alpha({\mathcal T}'_k)=-\sum_{T\in {\mathcal T}'_k}(-1)^{d^\circ(T)}\ 
Z_\alpha(T)$$
and 
$$\tilde Z_\alpha({\mathcal T}'_k)=\sum_{T\in {\mathcal T}'_k}(-1)^{d(T)}\ 
\tilde Z_\alpha(T)$$
where ${\mathcal T}'_k={\mathcal T}_k\setminus\{r\}$
denotes the set of all $k-$regular trees
which are non-trivial. We have obviously
$$Z({\mathcal T}_k)=-X+\sum_{\alpha\in{\mathcal A}} Z_\alpha({\mathcal T}'_k)
\qquad \hbox{ and }\qquad 
\tilde Z({\mathcal T}_k)=-X+\sum_{\alpha\in{\mathcal A}} \tilde Z_\alpha(
{\mathcal T}'_k)\ .$$

\begin{prop} \label{spinequality} We have 
$$Z_\alpha({\mathcal T}'_k)=Y_\alpha\left(X-(M_1\ W)_\alpha\right)\cdots
\left(X-(M_k\ W)_\alpha\right)$$
and 
$$\tilde Z_\alpha({\mathcal T}'_k)=Y_\alpha\left(-X+(\tilde M_1\ \tilde W)_\alpha\right)\cdots
\left(-X+(\tilde M_k\ \tilde W)_\alpha\right)
$$
where $W$ (respectively $\tilde W$) is the column vector with coordinates
$Z_\alpha({\mathcal T}'_k)$ (respectively 
$\tilde Z_\alpha({\mathcal T}'_k)$) indexed by $\alpha\in{\mathcal A}$.
\end{prop}

For a fixed natural integer $k\geq 1$ and $k$ (weight-)matrices
$M_1,\dots,M_k$ we have now:

\begin{cor} \label{g=Z} We have $g_\alpha=Z_\alpha({\mathcal T}'_k),\ 
\tilde g_\alpha=\tilde Z_\alpha({\mathcal T}'_k)$ for all 
$\alpha\in{\mathcal A}$ and 
$g=Z({\mathcal T}_k),\ 
\tilde g=\tilde Z({\mathcal T}_k)$
where $g_\alpha,\ \tilde g_\alpha,\ g,\ \tilde g$ are as in Theorem
\ref{main} and where $Z_\alpha({\mathcal T}'_k),\ \tilde 
Z_\alpha({\mathcal T}'_k),\ Z({\mathcal T}_k),\ \tilde Z({\mathcal T}_k)$
are as above.
\end{cor}

{\bf Proof of Proposition \ref{spinequality}} 
Given a non-trivial $k-$regular tree $T\not=\{r\}$ with root $r$,
we denote by $T_1,\dots,T_k$ the principal subtrees associated to 
the $k$ linearly ordered sons $s_1<s_2 \dots<s_k$ of $r$.

The generating function  
$$Z_\alpha({\mathcal T}'_k)=-\sum_{T\in {\mathcal T}'_k}
(-1)^{d^\circ(T)}\ Z_\alpha(T)=\sum_{T\in {\mathcal T}'_k}
(-1)^{\sum_jd^\circ(T_j)}\ Z_\alpha(T)$$
(with the last equality following from
$d^\circ(T)=1+\sum_{j=1}^k d^\circ(T_j)$) factorizes now as

$$\begin{array}{rcl}
\displaystyle Z_\alpha({\mathcal T}'_k)&=&Y_\alpha
\Big(X+\sum_{\beta\in{\mathcal A}} M_1(\alpha,\beta)\sum_{T\in {\mathcal T}'_k}
(-1)^{d^\circ(T)}Z_\beta(T)\Big)\cdot\\
&&\quad\cdot\Big(X+\sum_{\beta\in{\mathcal A}} M_2(\alpha,\beta)
\sum_{T\in {\mathcal T}'_k}(-1)^{d^\circ(T)}Z_\beta(T)\Big)\cdots\\
&&\quad\cdots\ \Big(X+\sum_{\beta\in{\mathcal A}} M_k(\alpha,\beta)
\sum_{T\in {\mathcal T}'_k}(-1)^{d^\circ(T)}Z_\beta(T)\Big)\\
\displaystyle &=&Y_\alpha
\Big(X-\sum_{\beta\in{\mathcal A}} M_1(\alpha,\beta)\ 
Z_\beta({\mathcal T}'_k\Big)\cdots
\Big(X-\sum_{\beta\in{\mathcal A}} M_k(\alpha,\beta)\ 
Z_\beta({\mathcal T}'_k)\Big)\ .\end{array}$$
Indeed, neglecting the trivial term $Y_\alpha$, 
the partition function $Z_\alpha(T)$ of a tree $T\in{\mathcal T}'_k$ 
with prescribed spin
$\varphi(r)=\alpha$ on its root decomposes into $k$ obvious factors
corresponding to the $k$ principal subtrees of $T$: A 
leafless edge $e\in E'$ with $\omega(e)\in T_j$ yields a
contribution to the $j-$th factor of $Z_\alpha({\mathcal T}'_k)$. Summing over
all $k-$regular trees $T\in {\mathcal T}'_k$ yields the above expression. The
term $X$ in the $j-$th factor corresponds to trees $T\in {\mathcal T}'_k$ 
whose $j-$th principal subtree $T_j$ is trivial.
This proves the first equality of Proposition \ref{spinequality}.

For 
$$\tilde Z_\alpha({\mathcal T}'_k)=\sum_{T\in {\mathcal T}'_k}
(-t)^{d(T)}\ \tilde Z_\alpha(T)=\sum_{T\in {\mathcal T}'_k}
(-t)^{\sum_j d(T_j)}\ \tilde Z_\alpha(T)$$
we get the analogous factorization
$$\begin{array}{rcl}
\displaystyle \tilde Z_\alpha({\mathcal T}'_k)&=&
Y_\alpha\Big(-X+\sum_{\beta\in{\mathcal A}} 
\tilde M_1(\alpha,\beta)\sum_{T\in {\mathcal T}'}(-1)^{d(T)}\tilde Z_\beta(T)
\Big)\cdot\\
&&\quad \cdot\Big(-X+\sum_{\beta\in{\mathcal A}} 
\tilde M_2(\alpha,\beta)\sum_{T\in {\mathcal T}'}(-1)^{d(T)}\tilde Z_\beta(T)
\Big)\cdots\\
&&\quad \cdots\Big(-X+\sum_{\beta\in{\mathcal A}} 
\tilde M_k(\alpha,\beta)\sum_{T\in {\mathcal T}'}(-1)^{d(T)}\tilde Z_\beta(T)
\Big)\\
\displaystyle &=&Y_\alpha\Big(-X+\sum_{\beta\in{\mathcal A}} \tilde 
M_1(\alpha,\beta)\ \tilde Z_\beta({\mathcal T}'_k)\Big)\cdots
\Big(-X+\sum_{\beta\in{\mathcal A}} \tilde 
M_k(\alpha,\beta)\ \tilde Z_\beta({\mathcal T}_k)\Big)
\end{array}
$$
which proves the second equality of Proposition \ref{spinequality}.
\hfill$\Box$

The proof of Corollary \ref{g=Z} is immediate: corresponding series
are recursively defined by the same formulae and initial data.

\section{Proof of Theorem \ref{main}}

By Corollary \ref{g=Z}, 
the formal power series $g=g(X)=-X+\sum_{\alpha\in{\mathcal A}}
g_\alpha(X)$ and $\tilde g=\tilde g(X)=-X+\sum_{\alpha\in{\mathcal A}}
\tilde g_\alpha(X)$ involved
in Theorem \ref{main} are suitably signed generating series 
for the partition
functions of complementary spin-models defined on $k-$regular
trees. In order to prove Theorem \ref{main} we define spin-models
on a set of combinatorial objects which we call grafted trees. 
A suitably signed
generating series of the partition functions for these spin-models
equals $g\circ_X\tilde g$ and a direct computation establishes
Theorem \ref{main}.

\subsection{A spin model on grafted trees}

A {\it grafted tree} is given by
$$(A;B_1,B_2,\dots,B_{d(A)})$$
where $A$ is a tree with $d(A)$ leaves
and where $B_1,\dots,B_{d(A)}$ is a sequence of $d(A)$ trees.
A grafted tree is {\it $k$-regular} if it involves only
$k-$regular trees. 

The {\it skeleton} of a grafted tree $(A;B_1,\dots,B_{d(A)})$ is the 
tree obtained by grafting (gluing) the root of
$B_1$ to the smallest (leftmost) leaf of $A$, by grafting the root of $B_2$ to
the second-smallest leaf of $A$ etc. The skeleton of $(A;B_1,\dots,B_{d(A)})$
has thus $\sum_{j=1}^{d(A)}d(B_j)$ leaves.

The data of a $k$-regular grafted tree $(A;B_1,\dots,B_{d(A)})$ and 
$k$ weight-matrices $M_1,\dots,M_k$ indexed by 
${\mathcal A}\times {\mathcal A}$ defines a
\lq\lq compositional'' spin model as follows: 
Compute first the ordinary partition function $Z(A)$ as defined
previously. Compute also the partition functions 
$\tilde Z(B_1),\dots,\tilde Z(B_{d(A)})$ with respect to the complementary
weight matrices $\tilde M_1=J-M_1,\dots,\tilde M_k=J-M_k$ (where
we denote by $J$ the all $1$ matrix indexed by 
${\mathcal A}\times {\mathcal A}$).
Replace now the $i-$th (occurence of the) letter $X$ in every monomial
of $Z(A)$ by $\tilde Z(B_i)$. The resulting formal power series
(or non-commutative polynomial for a finite
alphabet $\mathcal A$)
is by definition $Z(A;B_1,\dots,B_{d(A)})$. It is of degree
$\sum_{i=1}^{d(A)}d(B_i)$ in $X$.

In particular, we have 
$$Z(A;B_1,\dots,B_{d(A)})=X^{-d(V_A)}Z(A)\prod_{j=1}^{d(A)}\tilde Z(B_j)$$
(with $d(V_A)$ denoting the number of leaves in $A$)
in a commutative setting.

Let $g(X)=-X+\sum_{\alpha\in{\mathcal A}}\ g_\alpha(X),\ \tilde g(X)
=-X+\sum_{\alpha\in{\mathcal A}}\
\tilde g_\alpha(X)$ be the generating series involved in Theorem
\ref{main} associated to $k$ matrices $M_1,\dots,M_k$.

\begin{prop} \label{mainprop} We have
$$g\circ_X\tilde g(X)=-\sum_{(A;B_1,\dots,B_{d(A)})\in {\mathcal G}_k}
(-1)^{d^\circ(A)+\sum_{j=1}^{d(A)}d(B_j)}\ Z(A;B_1,\dots,B_{d(A)})$$
where ${\mathcal G}_k$ denotes the set of all $k-$regular grafted trees.
\end{prop}

{\bf Proof.} This follows at once from Corollary \ref{g=Z} and the 
definition of spin-models on grafted trees.\hfill$\Box$

Given a $k-$regular tree $T\in{\mathcal T}_k$
we denote by ${\mathcal S}(T)$ the set of all $k-$regular grafted trees
with skeleton $T$. Elements of ${\mathcal S}(T)$ are in bijection with 
$k-$regular rooted subtrees of $T$ containing the root $r\in T$: 
The subtree $A\subset T$ of
a $k-$regular grafted tree $(A;B_1,\dots,B_{d(A)})\in
{\mathcal S}(T)$ with skeleton $T$ clearly contains $r$.
Since the subtree $B_j\subset T$ is the maximal subtree of $T$
rooted in the $j-$th leaf $v_j$ of $A$, the pair of trees $A\subset T$ 
defines the grafted tree $(A;B_1,\dots,B_{d(A)})$ completely. 

We denote by ${\mathcal S}'(T)={\mathcal S}\setminus
\{(\{r\};T)\}$ the set of all 
$k-$regular grafted trees $(A;B_1,\dots,B_{d(A)})$ with skeleton $T$ and
non-trivial $A\not=\{r\}$.

For $\alpha\in{\mathcal A}$ and 
$(A;B_1,\dots,B_{d(A)})\in {\mathcal S}'(T)$ we define
$Z_\alpha(A;B_1,\dots,B_{d(A)})$ in the obvious way
by replacing the $i-$th occurence of the letter $X$ in 
$Z_\alpha(A)$ with $\tilde Z(B_i)$ for $i=1,\dots,d(A)$.

\begin{prop} \label{equalzero} We have for a non-trivial $k-$regular tree 
$T\in{\mathcal T}'_k$ and for all $\alpha\in {\mathcal A}$
$$\sum_{(A;B_1,\dots,B_{d(A)})\in{\mathcal S}'(T)}(-1)^{d^\circ(A)}
Z_\alpha(A;B_1,\dots,B_{d(A)})=-\tilde Z_\alpha(T)\ .$$
\end{prop}

{\bf Proof.} The proof is by
induction on the number $d(T)$ of leaves in $T$. If all $k$ principal 
subtrees associated to the linearly ordered sons
$s_1<\dots<s_k$ of the root $r\in T$ are trivial, we have
${\mathcal S}'(T)=(T;\{s_1\},\dots,\{s_k\})$ and
$$\begin{array}{cl}
\displaystyle &\sum_{(A;B_1,\dots,B_{d(A)})\in{\mathcal S}'(T)}(-1)^{d^\circ(A)}
Z_\alpha(A;B_1,\dots,B_{d(A)})\\
\displaystyle =&
(-1)^{d^\circ(T)}Z_\alpha(T;\{s_1\},\dots,\{s_k\})=-\ Y_\alpha\ .\end{array}$$
Since we have $\tilde Z_\alpha(T)=Y_\alpha$, Proposition 
\ref{equalzero} follows in this case.

Otherwise, we have
$$\sum_{(A;B_1,\dots,B_{d(A)})\in
{\mathcal S}'(T)}(-1)^{d^\circ(A)}Z_\alpha(A;B_1,\dots,
B_{d(A)})=-Y_\alpha C_1 C_2\cdots C_k$$
where $C_i=X$ if the principal subtree $T_i=\{s_i\}$ is trivial and 
$$C_i=\tilde Z(T_i)+\sum_{\beta\in{\mathcal A}}M_i(\alpha,\beta)
\sum_{(A;B_1,\dots,B_{d(A)})\in {\mathcal S}'(T_i)}
Z_\beta(A;B_1,\dots,B_{d(A)})$$
if $T_i\not=\{s_i\}$. The contribution of a non-trivial principal
subtree $T_i\not=\{s_i\}$ is thus by induction
$$\begin{array}{rcl}
C_i&=&\tilde Z(T_i)-\sum_{\beta\in{\mathcal A}}M_i(\alpha,\beta)
\tilde Z_\beta(T_i)\\
&=&\sum_{\beta\in{\mathcal A}} \left(\tilde Z_\beta(T_i)-M_i(\alpha,\beta)
\tilde Z_\beta(T_i)\right)\\
&=&\sum_{\beta\in{\mathcal A}} \left(1- M_i(\alpha,\beta)\right)
\tilde Z_\beta(T_i)\\
&=&\sum_{\beta\in{\mathcal A}} \tilde M_i(\alpha,\beta)\tilde Z_\beta(T_i)
\end{array}$$
and
$$-Y_\alpha C_1C_2\cdots C_k=-\tilde Z_\alpha(T)\ .\qquad \Box$$

\begin{cor} \label{eqcor}
We have for a non-trivial $k-$regular tree $T\in{\mathcal T}'_k$
$$\sum_{(A;B_1,\dots,B_{d(A)})\in{\mathcal S}(T)}
(-1)^{d^\circ(A)}Z(A;B_1,\dots,B_{d(A)})=0\ .$$
\end{cor}

{\bf Proof.} Sum the equality of Proposition \ref{equalzero}
over $\alpha\in{\mathcal A}$ and use 
$$Z(\{r\};T)=\tilde Z(T)=\sum_{\alpha\in{\mathcal A}} 
\tilde Z_\alpha(T)\ .\qquad \Box$$

{\bf Proof of Theorem \ref{main}.} By Proposition \ref{mainprop}
we have 
$$g\circ_X\tilde g=-\sum_{T\in{\mathcal T}_k}(-1)^{d(T)}
\sum_{(A;B_1,\dots,B_{d(A)})\in{\mathcal S}(T)}(-1)^{d^\circ(A)}
Z(A;B_1,\dots,B_{d(A)})\ .$$

It follows from Corollary \ref{eqcor} that only the trivial tree
$T=\{r\}$ contributes to the right hand side.
The contribution of $T=\{r\}$ amounts to 
$$-(-1)^{d(\{r\})}Z(\{r\};\{r\})=-(-1)\cdot X=X$$
which proves Theorem \ref{main}.\hfill$\Box$


\section{Further generalizations}

We give first a generalization of Theorem \ref{main} associated
to not necessarily regular trees. The special case of this generalization
where all weight-matrices have coefficients in $\{0,1\}$ 
appears already in \cite{P}.

We state then an even more general version involving trees with vertices
labelled by a set ${\mathcal T}$ of types and formal power-series 
in $R[[(X_\tau)_{\tau\in{\mathcal T}},(Y_\alpha)_{\alpha\in{\mathcal A}}]]$.

\subsection{Trees with vertices of arbitrary degrees}

Consider a subset ${\mathcal K}\subset {\mathbb N}_{>0}$ of strictly 
positive integers and a partition ${\mathcal A}=\cup_{k\in{\mathcal K}}
{\mathcal A}_k$ of the alphabet $\mathcal A$ into non-empty parts indexed
by ${\mathcal K}$. 
For each $k\in 
{\mathcal K}$ choose $k$ matrices $M_1^k,\ M_2^k,\dots, M_k^k$
with rows indexed by ${\mathcal A}_k$ and columns indexed by
${\mathcal A}$. For $i=1,\dots,k$, set $\tilde M_i^k=J-M_i^k$ where 
$J$ is the all $1$ matrix with indices in ${\mathcal A}_k\times {\mathcal A}$.

For $\alpha\in{\mathcal A}_k$ 
consider the uniquely defined
(non-commutative) formal power series $g_\alpha,\ \tilde g_\alpha$
consisting only of monomials of degree at least $1$ in $X$, involving
at least one variable of $(Y_\beta)_{\beta\in {\mathcal A}}$ and such that
$$g_\alpha=Y_\alpha(X-(M_1^k\ V)_\alpha)(X-(M_2^k\ V)_\alpha)\cdots
(X-(M_k^k\ V)_\alpha)$$
and 
$$\tilde g_\alpha=Y_\alpha(-X+(\tilde M_1^k\ \tilde V)_\alpha)\cdots
(-X+(\tilde M_k^k\ \tilde V)_\alpha)$$
where $V$, respectively $\tilde V$, is the column vector with coordinates
$g_\beta$, respectively $\tilde g_\beta$, for $\beta\in{\mathcal A}$. 

Set $g=-X+\sum_{\alpha\in {\mathcal A}} g_\alpha$ and 
$\tilde g=-X+\sum_{\alpha\in {\mathcal A}} \tilde g_\alpha$.

\begin{thm} \label{nonregular}
We have 
$$g\circ_X \tilde g=\tilde g\circ_X g=X\ .$$
\end{thm}

The proof is an adaption of the proof of Theorem \ref{main}: 
Consider trees whose internal vertices have degrees in ${\mathcal K}$
and colour internal vertices of degree $k\in{\mathcal K}$ by elements
in ${\mathcal A}_k$.

{\bf Example (Inversion of power series).} 
In this example, we work over a commutative
ring $R$ and set $Y_\alpha=1$. Choose an integer $l\geq 1$ and 
set ${\mathcal K}=\{2,3,4,\dots\}$ 
(this ensures existence of all relevant 
power-series). We consider trees with internal 
vertices of degree $\geq 2$ having spins in a finite alphabet 
${\mathcal A}$ containing $l$ elements. 

For $1\leq j\leq k,\ 2\leq k$ choose  elements $\beta_{k,j}\in R$ and
consider the diagonal
weight function $$w_e(\varphi(\alpha(e)),\varphi(\omega(e))=
\beta_{k,j}\ \delta_{\varphi(\alpha(e)),\varphi(
\omega(e))}=\left\{\begin{array}{ll} \beta_{k,j}&\hbox{if }
\varphi(\alpha(e))=\varphi(
\omega(e))\\
0&\hbox{otherwise}\end{array}\right.$$
where the leafless edge $e$ joins a vertex $\alpha(e)$ of degree $k$ to 
its $j-$th son $\omega(e)$.

Introduce the formal power series $g_*,\ \tilde g_*\in R[[X]]$ 
without constant term which satisfy
$$\begin{array} {lcl}
g_*&=&\sum_{k=2}^\infty \ \prod_{j=1}^k\left(X-\beta_{k,j}\ g_*\right)\\
\tilde g_*&=&\sum_{k=2}^\infty \ \prod_{j=1}^k\left(-X+(l-\beta_{k,j})
\ \tilde g_*\right)\ .\end{array}$$

Set $g(X)=-X+l\ g_*$ and $\tilde g(X)=-X+l\ \tilde g_*$.
Theorem \ref{nonregular}, applied to the present situation shows that
we have
$$g(\tilde g(X))=\tilde g(g(X))=X\ .$$
Since this formula holds for any $l\in{\mathbb N}$, it extends to 
an arbitrary value of $l\in R$ which can thus be considered as
a parameter. 

This formula provides a mean (other
than the celebrated Lagrange inversion formula)
for computing the compositional
inverse of a formal power series $h(X)=\sum_{k=1}^\infty \gamma_k\ X^k$
with $\gamma_1\in R^*$ invertible: Set $f=h(-X/\gamma_1)=
-X+\sum_{k=2}^\infty \beta_k\ t^k$. If $\beta_2\not=0$ set 
$l=\beta_2$ and choose
constants $\beta_{k,j}$ such that $g(X)=f(X)$ with $g(X)$ defined as above.

Each such choice can be used to compute the compositional
inverse $\tilde g=f^{-1}$ of $f$. The compositional inverse of the 
initial series
$h(X)$ is then given by $-\gamma_1\ \tilde g(X)$. The case $\beta_2=0$
can be handled similarly by allowing the trees to have internal vertices of
degree $1$.

\subsection{Vertices of different types}

This still more general version of Theorem \ref{main} is perhaps
easier to formulate in a combinatorial way:

Consider a set ${\mathcal T}$ of {\it vertex-types} indexing the 
non-empty parts of a partition $
{\mathcal A}=\cup_{\tau\in{\mathcal T}} {\mathcal A}_\tau$. 
Consider also a set ${\mathcal E}$ of
{\it edge-types} together with weight-functions
$$w_\epsilon:{\mathcal A}\times{\mathcal A}\longrightarrow R$$
indexed by edge-types $\epsilon\in {\mathcal E}$.

A {\it labelled} tree is a finite rooted planar tree $L$ with leaves
labelled by elements of ${\mathcal T}$, internal vertices labelled by 
pairs $(\tau,\alpha\in {\mathcal A}_\tau)\in {\mathcal T}\times {\mathcal A}$
and edges labelled by elements of ${\mathcal E}$. A vertex $v\in L$ is of type
$\tau\in{\mathcal T}$ if it is either a leaf labelled $\tau$ or
an internal vertex labelled $(\tau,\alpha\in{\mathcal A}_\tau)$.
A labelled tree $L$ is of type $\tau\in{\mathcal T}$ if its root vertex is 
of type $\tau$.

The {\it energy} $f(L)\in R[(X_\tau)_{\tau\in {\mathcal T}},
(Y_\alpha)_{\alpha\in{\mathcal A}}]$
of a labelled tree $L$ of type $\tau\in{\mathcal T}$
is defined in the obvious way:
$f(L)=X_\tau$ if $L$ is trivial (reduced to its labelled root)
and
$$f(L)=Y_\alpha C_1\cdots C_k$$
otherwise where $(\tau,\alpha\in{\mathcal A}\tau)$ is
the label of the root $r\in L$ and where
 $C_i\in R[(X_\tau)_{\tau\in {\mathcal T}},
(Y_\alpha)_{\alpha\in{\mathcal A}}]$ 
is associated to the $i-$th labelled principal 
subtree $L(s_i)$ of $L$ as follows: $C_i=X_{\tau_i}$ if 
$L(s_i)$ is trivial of type $\tau_i$ and
$$C_i=
w_{\epsilon_i}(\alpha,\alpha_i)f(L_i)$$
where $\epsilon_i$ is the label of the edge joining the root $r$ to
its $i-$th son $s_i$ which is labelled 
$(\tau_i,\alpha_i\in{\mathcal A}_{\tau_i})$.

Define complementary weight-functions by $\tilde w_\epsilon(\alpha,\beta)=
1-w_\epsilon(\alpha,\beta)$ and compute the {\it complementary energy}
$\tilde f(L)$ of $L$ as above using the complementary weight-functions.

Call two labelled trees $L_1,L_2$ {\it edge-equivalent} if they differ only 
on their edge-labels (but share the same underlying tree-structure
and vertex-labels).

Let ${\mathcal L}$ be a subset of the set of all labelled (finite
rooted planar) trees. We denote by
${\mathcal L}_\tau\subset {\mathcal L}$ the subset of all labelled trees 
of type $\tau$ in ${\mathcal L}$ and assume that the set
${\mathcal L}$ satisfies the following three conditions:

\ \ (i) If $L\in{\mathcal L}$ and $v$ is a leaf of ${\mathcal L}$,
then $L(v)\in{\mathcal L}$ where $L(v)$ denotes the labelled
subtree defined in the obvious way by considering the 
maximal subtree of $L$ rooted at $v$. 

\ \ (ii) If $L'\in{\mathcal L}_\tau$ and $v$ is a leaf of type $\tau$
in $L\in{\mathcal L}$ then the labelled tree obtained in the obvious 
way by gluing $L'$ onto the leaf $v\in L$ is again in $\mathcal L$.

\ \ (iii) All equivalence classes of the edge-equivalence relation are
finite.

Otherwise stated, all maximal labelled subtrees of an element in 
${\mathcal L}$ are in ${\mathcal L}$ by (i) and ${\mathcal L}$ is
\lq\lq closed under composition'' by (ii).
Condition (iii) is a finiteness condition (which can perhaps be
slightly weakened or replaced by a similar statement) 
ensuring the existence
of the generating series $Z_\tau$ and $\tilde  Z_\tau$ defined below.

For such a set ${\mathcal L}$ we define
$$Z_\tau=-\sum_{L\in{\mathcal L}}(-1)^{d^\circ(L)} f(L)$$
and 
$$\tilde  Z_\tau=\sum_{L\in{\mathcal L}}(-1)^{d(L)} \tilde f(L)$$
where $d^\circ(L)$ denotes the number of internal vertices of a
labelled tree $L$ and where $d(L)$ denotes the number of leaves in
$L$.

\begin{thm} We have
$$Z_\tau\circ_{X_{\mathcal L}}(\tilde Z_\sigma)_{\sigma\in{\mathcal T}}=
\tilde Z_\tau\circ_{X_{\mathcal L}}(Z_\sigma)_{\sigma\in{\mathcal T}}
=X_\tau$$
where the notation $Z_\tau\circ_{X_{\mathcal L}}
(\tilde Z_\sigma)_{\sigma\in{\mathcal T}}$ means that every occurence of
$X_\tau$ in $Z_\tau$ is replaced by the generating series
$\tilde Z_\tau$.
\end{thm}

This result can be used for the formal inversion of power-series in
several variables. 

{\bf Sketch of Proof.} Define grafted labelled trees in the obvious
way and check that Proposition \ref{equalzero} remains valid in the 
present context.



\section{Morphisms of rooted trees into posets}

This section is a digression discussing a generalization of 
grafted trees.

The vertex set $V$ of a rooted (not necessarily finite or planar) 
tree $T$ with oriented edges $E$ can be considered as a poset 
(partially ordered set)
by considering the order relation induced by $\omega(e)>\alpha(e)$ for 
$e\in E$. 

\begin{rem} Considering a rooted tree as a poset, one might wonder how
many unordered pairs of comparable vertices are contained in 
($k-$regular) rooted planar trees having $n$ vertices. Let $a_n$
(respectively $a_{n,k}$) denote this number. We have then
$$\sum_{n=1}^\infty a_n\ t^n=\frac{\partial y}{\partial u}(t,1)$$
where $y(t,u)=t+(\hbox{ terms of higher order in }t\ )$
satisfies the functional equation
$$y(t,u)=t+t\frac{y(tu,u)}{1-y(tu,u)}\ .$$
The corresponding sequence $a_2,a_3,a_4,\dots$ starts as 
$$1,5,22,93,386,1586,6476,26333,106762,431910,\dots\hbox{ cf. A346 of 
\cite{EIS}}.$$

Similarly, for the numbers $a_{n,k}$ associated to $k-$regular 
trees, we have
$$\sum_{n=1}^\infty a_{n,k}\ t^n=\frac{\partial y}{\partial u}(t,1)$$
with $y(t,u)=t+(\hbox{ terms of higher order in }t\ )$
satisfying
$$y(t,u)=t+t\ (y(tu,u))^k\ .$$
For $k=2$ we get the sequence $\frac{1}{2}(a_{3,2},a_{5,2},a_{7,2},\dots)$
starting as 
$$1,6,29,130,562,2380,9949,41226,169766,\dots ,\hbox{ cf. A8549 of 
\cite{EIS}}$$
and for $k=3$, the sequence $\frac{1}{3}(a_{4,3},a_{7,3},a_{10,3},\dots)$
starting as
$$1,9,69,502,3564,24960,173325,1196748,\dots ,\hbox{ cf. A75045 of 
\cite{EIS}.}$$
\end{rem}

Given a second poset $P$, a {\it morphism} from $T$ to $P$ is an 
application $\mu:V\longrightarrow P$ such that
$\mu(\omega(e))\geq \mu(\alpha(e))$ for every edge $e$ of $T$. 

Denote by $\{1,2\}$ (with $1<2$) the obvious totally ordered poset.
Call a morphism $\mu:T\longrightarrow\{1,2\}$ {\it restricted} if 
$\mu^{-1}(2)$ contains all leaves of $T$.

\begin{prop} \label{graftedmorphism}
If $T$ is a rooted (planar) tree, then the set of
grafted trees $(A;B_1,\dots,B_{d(A)})$ with skeleton $T$ 
corresponds bijectively to the set of restricted morphisms 
$\mu:T\longrightarrow \{1,2\}$. 
\end{prop}

{\bf Proof.} Given a restricted morphism $\mu:T\longrightarrow\{1,2\}$ 
we set $A=\{\mu^{-1}(1)\cup\ \hbox{sons of }\mu^{-1}(1)\}$ where
the root $r\in T$ is by convention the son of the empty set.
The rooted trees $\{B_1,\dots,B_d(A)\}$ are the connected components of 
the maximal subforest with vertices $\mu^{-1}(2)$. 
It is obvious that $(A;B_1,\dots,B_{d(A)})$
is a grafted tree with skeleton $T$. Reciprocally, given a grafted
tree $(A;B_1,\dots,B_{d(A)})$ with skeleton $T$, we get a
morphism $\mu:T\longrightarrow \{1,2\}$ by setting $\mu(v)=2$ if 
$v\in B_j$ for some $j=1,\dots,d(A)$ and $\mu(v)=1$ otherwise. 
Since the forest
$\cup_{j=1}^{d(A)}B_j$ contains all leaves of $T$, the morphism $\mu$ 
is restricted.
\hfill$\Box$

Proposition \ref{graftedmorphism} suggests to generalize
the notion of grafted trees by considering sets of (suitable)
morphisms from rooted trees into posets.
We consider now a few special cases. Henceforth all
rooted trees and posets will be finite. The number of morphisms of a 
rooted tree $T$ having $\sharp(V)$ vertices into a poset 
$P$ with $\sharp(P)$ elements is bounded by $\sharp(P)^{\sharp(V)}$.

\subsection{Grafted trees}

In this subsection we indicate how to 
count the number of grafted trees with given skeleton $T$.
A slight modification counts
all morphisms from $T$ into $\{1,2\}$. Further generalizations
consist in counting all (suitable) morphisms from $T$ into the totally 
ordered set $\{1,\dots,n\}$.

Given a (finite rooted planar) tree $T$ we denote by $\gamma(T)$
the number of grafted trees with skeleton $T$. We denote by 
$\tilde \gamma(T)\geq \gamma(T)$ the number
of all morphisms $T\longrightarrow \{1,2\}$.
We have $\gamma(\{r\})=1,\ \tilde \gamma(\{r\})=2$ if $T=\{r\}$ is the trivial
tree reduced to its root. Remark that $\tilde\gamma(\tilde T)=
\gamma(T)$ where $\tilde T$ is obtained by removing all leaves from
the rooted tree $\gamma$ (where $\tilde\gamma(\emptyset)=1$ by convention). 
This remark generalizes easily to the analogous numbers enumerating 
(restricted) morphismes $T\longrightarrow\{1,\dots,m\}$.

\begin{prop} \label{graftcount} For $T$ a non-trivial
(finite rooted planar) tree we have
$$\gamma(T)=1+\prod_{j=1}^k\gamma(T_i)$$
and
$$\tilde \gamma(T)=1+\prod_{j=1}^k\tilde \gamma(T_i)$$
where $T_1,\dots,T_k$ are the principal subtrees defined
by the $k$ sons $s_1,\dots,s_k$ of the root $r\in T$.
\end{prop}

{\bf Proof.} We count all (restricted)
morphisms $\mu$ from $T$ into $\{1,2\}$.

Consider a morphism $\mu:V\longrightarrow \{1,2\}$. 
If $\mu(r)=2$ we have $\mu(v)=2$ for all $v\in V$ and there
is exactly one such morphism yielding a contribution of $1$ to 
$\gamma(T)$ and $\tilde \gamma(T)$. If $\mu(r)=1$, the 
restriction $\mu_i$ of $\mu$ to a principal subtrees $T_i$ is an
arbitrary (restricted) morphism from $T_i$ into $\{1,2\}$ and
the restrictions $\mu_1,\dots,\mu_k$ can be
arbitrary thus proving the formula.\hfill$\Box$

Proposition \ref{graftcount} leads to a fast algorithm 
for computing $\gamma(T)$ illustrated by the following example.

{\bf Example.} For the tree $T$ of Figure 2
we have
$$\begin{array}{l}\gamma(v_1)=\gamma(v_3)=\gamma(v_8)=2\\
\gamma(v_4)=1+\gamma(v_3)\cdot 1=3,\ \gamma(v_7)=1+1\cdot \gamma(v_8)=3\\
\gamma(v_2)=1+\gamma(v_1)\cdot\gamma(v_4)=7,\ \gamma(v_6)=1+1\cdot \gamma(v_7)
=4\\
\gamma(v_5)=1+\gamma(v_2)\cdot\gamma(v_6)=29\end{array}$$
and 
$$\begin{array}{l}\tilde \gamma(v_1)=\tilde \gamma(v_3)=\tilde \gamma(v_8)=1+
2\cdot 2=5\\
\gamma(v_4)=1+\gamma(v_3)\cdot 2=11,\ \gamma(v_7)=1+2\cdot \gamma(v_8)=11\\
\gamma(v_2)=1+\gamma(v_1)\cdot\gamma(v_4)=56,\ \gamma(v_6)=1+2\cdot \gamma(v_7)
=23\\
\gamma(v_5)=1+\gamma(v_2)\cdot\gamma(v_6)=1289\end{array}$$
which shows that the tree $T$ of Figure 2
is the skeleton of $\gamma(T)=\gamma(v_5)=29$
different grafted trees (or equivalently, that $T$ contains 
$29$ different $2-$regular subtrees sharing the root $r$ with $T$)
and has $1289$ distinct morphisms into the totally
ordered poset $\{1,2\}$.

Proposition \ref{graftcount} shows also that the generating function 
$$y(t)=\sum_{T\in{\mathcal T}_2}\gamma(T)\ t^{d(T)}$$
counting the number $\gamma(k)=[t^k]y$ 
of $2-$regular grafted trees whose skelettons have $k$ leaves satisfies
$$\begin{array}{rcl}
\displaystyle y&=&\frac{1-\sqrt{1-4t}}{2}+y^2=\frac{1-\sqrt{-1+2\sqrt{1-4t}}}
{2}\\
\displaystyle &=&t+2t^2+6t^3+21t^4+80t^5+322t^6+1348t^7+5814t^8+\dots\ .
\end{array}$$

The similar generating function
$\tilde y(t)=\sum_{T\in{\mathcal T}_2}\tilde \gamma(T)\ t^{d(T)}$
(which counts all morphisms of binary rooted trees into $\{1,2\}$)
is given by
$$\begin{array}{rcl}
\displaystyle \tilde y&=&t+\frac{1-\sqrt{1-4t}}{2}+y^2=\frac{1-
\sqrt{-1-4t+2\sqrt{1-4t}}}{2}\\
\displaystyle &=&2t+5t^2+22t^3+118t^4+706t^5+4530t^6+\dots\end{array}
\ .$$

More generally, consider a fixed natural integer $k\geq 2$. We call a morphism 
$\mu:T\longrightarrow \{1,\dots, k\}$ {\it restricted} if 
$\mu^{-1}$ contains all leaves of the rooted tree $T$. We denote
by $y_m(t)$ the generating function counting the number of (non-isomorphic) 
restricted morphisms from a $k-$regular tree having $n$ leaves into
the completely ordered finite set $\{1,\dots,m\}$. 
Similarly, we denote by $\tilde y_m(t)$
the analogous generating function 
counting all morphisms without restriction on the leaves.

\begin{prop} We have
$$y_1(t)=\tilde y_1(t)=t+\left(y_1(t)\right)^k$$
$$y_m(t)=t+\sum_{h=1}^m \left(y_h(t)\right)^k$$
$$\tilde y_m(t)=m t+\sum_{h=1}^m \left(\tilde y_h(t)\right)^k$$
\end{prop}

{\bf Sketch of proof.} The trivial tree $T$ yields a contribution of $t$
to $y_k$ and of $mt$ to $\tilde y(t)$. For a morphism
$\mu$ from a nontrivial tree into $\{1,\dots,m,\}$ 
we consider the image $\mu(r)$ of its root and the remaining 
possibilities of the induced restrictions
$\mu_1,\dots,\mu_k$ on the principal subtrees 
$T_1,\dots,T_k$ of $T$.\hfill$\Box$

For arbitrary finite planar non-empty trees we consider
the generating function counting the number of 
different restricted morphisms $\mu$ from a tree having $n$ vertices into
the completely ordered finite set $\{1,\dots,m\}$. The generating function 
$\tilde y_m(t)$ counts all such morphisms into $\{1,\dots,m\}$. One can then
prove the following result.

\begin{prop} We have
$$y_1(t)=\tilde y_1(t)=t+\frac{ty_1}{1-y_1}$$
$$y_m(t)=t+t\ \sum_{h=1}^m \frac{y_h}{1-y_h}$$
$$\tilde y_m(t)=m t+t\ \sum_{h=1}^m \frac{\tilde y_h}{1-\tilde y_h}$$
\end{prop}

The first instances are
$y_1=\tilde y_1=\frac{1-\sqrt{1-4t}}{2}$
(defining the Catalan numbers, cf. A108 of \cite{EIS})
$$\begin{array} {lcl}
\displaystyle 
y_2&=&\displaystyle 
\frac{3-2t-\sqrt{1-4t}-\sqrt{2-16t+4t^2+(2+4t)\sqrt{1-4t}}}{4}\\
&=&\displaystyle 
t+2t^2+5t^3+15t^4+50t^5+178t^6+663t^7+2553t^8+\dots\end{array}
$$
which is sequence  A7853 of \cite{EIS} and
$$\begin{array} {lcl}
\displaystyle \tilde y_2&=&\displaystyle
\frac{3-\sqrt{1-4t}-\sqrt{2-20t+2\sqrt{1-4t}}}{4}\\
&=&\displaystyle 
2t+3t^2+9t^3+34t^4+145t^5+667t^6+3231t^7+16247t^8+\dots\end{array}
$$

Finally, let us mention the well-known fact that the function 
$p(m)$ counting all morphisms
from a fixed finite poset $E$ into the totally ordered poset $\{1,\dots,m\}$
is a polynomial of degree $\sharp (E)$. Indeed
$$p(m)=\sum_{k=1} \alpha_k {m\choose k}$$
where $\alpha_k$ denotes the number of surjective morphismes of
$E$ into $\{1,\dots,k\}$. One can thus also consider 
generating functions associated to such polynomials. 
The number of restricted
morphisms from a fixed rooted tree into $\{1,\dots,m\}$ is
of course also a polynomial function. Its degree in $m$ is 
the number of interior leaves in $T$.

\subsection{Surjective morphisms} Given a finite rooted tree $T$ having
$n$ vertices, the
number $\sigma(T)$ of surjective morphisms from $T$ into 
$\{1,\dots,n\}$ can be recursively computed by remarking that
$$\sigma(T)=(n-1)!\ \prod_{j=1}^k \frac{\sigma(T_j)}{n_j!}$$
where $T_1,\dots,T_k$ are the principal (rooted) subtrees of $T$ 
having $n_1,\dots,n_k$ vertices.
Denoting by ${n-1\choose n_1,\dots,n_k}=\frac{(n-1)!}{n_1!\cdots n_k!}$
the obvious multinomial coefficient (where $n-1=\sum n_i$)
we get for our favorite (binary) tree $T$ of Figure $2$ the following 
numbers:
$\sigma(v)=1$ if $v$ is a leaf and
$$\sigma(v_1)=\sigma(v_3)=\sigma(v_8)={2\choose 1,1}=2$$
$$\sigma(v_4)=\sigma(v_7)={4\choose 1,3}\ 2=8$$
$$\sigma(v_2)={8\choose 3,5}\ 2\cdot 8=896,\ 
\sigma(v_6)={6\choose 1,5}\ 1\cdot 8=48$$
$$\sigma(T)=\sigma(v_5)={16\choose 9,7}\ 896\cdot 48=492011520$$ 

The generating function 
$$y(t)=\sum_{k=1}^\infty \alpha_kt^k$$
encoding the number $\alpha_n$ of surjective morphisms into
$\{1,\dots,n\}$ from all rooted binary planar trees
on $n$ vertices satisfies
$$\alpha_n=(n-1)!\sum_{k=1}^{n-2}\frac{\alpha_k\ \alpha_{n-1-k}}{k!\ (n-1-k)!}$$
(cf. sequence A182 of \cite{EIS}).
Otherwise stated, the exponential generating function 
$$z(t)=\sum_{k=1}^\infty \alpha_k\ \frac{t^k}{k!}$$
satisfies $z'=z^2-1$ thus proving that $z(t)=\hbox{tanh}(t)$.

Similarly, considering the exponential generating function 
$$z(t)=\sum_{k=1}^\infty \beta_k\ \frac{t^k}{k!}$$
enumerating the number $\beta_n$ of all surjective homomorphisms from  
rooted planar trees having $n$ vertices into $\{1,\dots,n\}$
we have 
$$z'=\frac{z}{1-z}+1=\frac{1}{1-z}$$
which implies $z(t)=\frac{1-\sqrt{1-4t}}{2}$ (cf. sequence 
A108 of \cite{EIS}) and shows that
$\beta_n=n!{2(n-1)\choose n-1}/n=\frac{(2n-2)!}{(n-1)!}$. 




\section{Loday's example (i)}

The aim of this section is a partial analysis of example (i) in 
\cite{L}. The framework is somewhat simpler as in the previous sections:
We work over the commutative ground field ${\mathbf C}$ of complex numbers.
The alphabet ${\mathcal A}$ consists of nine elements and equals
$${\mathcal A}=\{\circ,N,{NW},W,{SW},S,{SE},E,{NE}\}$$
suggesting the graphical notations of \cite{L}. We set $X=t$ in order
to stick to \cite{L}.

Loday's example (i) 
corresponds to $k=2$. The two $9\times 9$ matrices
$M_1,M_2$ are given by 
$$\left(\begin{array}{ccccccccc}
1&0&0&0&0&0&0&0&0\\
1&1&1&0&0&0&0&0&0\\
1&1&1&1&0&0&0&0&0\\
1&0&1&1&0&0&0&0&0\\
1&1&1&1&1&0&0&0&0\\
1&1&0&0&0&0&0&0&0\\
1&1&1&1&0&1&0&0&0\\
1&1&0&1&0&0&0&0&0\\
1&1&1&1&0&0&0&0&0\end{array}\right)\ ,\ 
\left(\begin{array}{ccccccccc}
0&1&1&1&0&0&0&0&0\\
0&0&0&1&0&0&0&0&0\\
0&0&0&0&0&0&0&0&0\\
0&1&0&0&0&0&0&0&0\\
0&0&0&0&0&1&0&0&0\\
0&0&1&1&1&1&0&0&0\\
0&0&0&0&1&0&1&1&1\\
0&0&1&0&0&0&0&1&1\\
0&0&0&0&0&0&0&1&1\end{array}\right)$$
Writing $V=(g_\circ,g_N,g_{NW},g_W,g_{SW},g_S,g_{SE},g_E,g_{NE})^t$
we get the equations
$$\begin{array}{lcl}
g_\circ&=&(-t+g_\circ)(-t+g_N+g_{NW}+g_W)\\
g_N&=&(-t+g_\circ+g_N+g_{NW})(-t+g_W)\\
g_{NW}&=&(-t+g_\circ+g_N+g_{NW}+g_W)(-t)\\
g_W&=&(-t+g_\circ+g_{NW}+g_W)(-t+g_N)\\
g_{SW}&=&(-t+g_\circ+g_N+g_{NW}+g_W+g_{SW})(-t+g_S)\\
g_S&=&(-t+g_\circ+g_N)(-t+g_{NW}+g_W+g_{SW}+g_S)\\
g_{SE}&=&(-t+g_\circ+g_N+g_{NW}+g_W+g_S)\\
& &\qquad (-t+g_{SW}+g_{SE}+g_E+g_{NE})\\
g_E&=&(-t+g_\circ+g_N+g_W)(-t+g_{NW}+g_E+g_{NE})\\
g_{NE}&=&(-t+g_\circ+g_N+g_{NW}+g_W)(-t+g_E+g_{NE})
\end{array}$$
One sees easily that we have $$g_o=g_N=g_W\ .$$
Eliminating $g_N,g_W$ we get the simpler
equations:
$$\begin{array}{lcl}
g_\circ&=&(-t+g_\circ)(-t+2g_\circ+g_{NW})\\
g_{NW}&=&(-t+3g_\circ+g_{NW})(-t)\\
\hline
g_{SW}&=&(-t+3g_\circ+g_{NW}+g_{SW})(-t+g_S)\\
g_S&=&(-t+2g_\circ)(-t+g_\circ+g_{NW}+g_S+g_{SW})\\
\hline
g_E&=&(-t+3g_\circ)(-t+g_{NW}+g_E+g_{NE})\\
g_{NE}&=&(-t+3g_\circ+g_{NW})(-t+g_E+g_{NE})\\
\hline
g_{SE}&=&(-t+3g_\circ+g_{NW}+g_S)(-t+g_{SW}+g_{SE}+g_E+g_{NE})\\
\end{array}$$
where functions above a horizontal line are independent 
of functions below the line.

Computations (done with Maple) using Groebner-bases
show that $$y=g_\circ+g_N+g_{NW}+g_W+g_{SW}+g_S+g_{SE}+g_E+g_{NE}$$
satisfies the algebraic equation $P(y(t),t)=0$ where
$$P(y,t)=c_0+c_1\  y+c_2\ y^2+c_3\ y^3+c_4\ y^4$$
with
$$\begin{array}{lcl}c_0&=&
t\ \left( 288\,{t}^{31}-1008\,{t}^{30}-17696\,{t}^{29}+35124\,{t}^{28}+513042\,{t}^{27}\right.\\&&
-352654\,{t}^{26}-8834409\,{t}^{25}-2315100\,{t}^{24}+
94293622\,{t}^{23}\\&&
+92841847\,{t}^{22}-608228325\,{t}^{21}-1031578684\,
{t}^{20}\\&&
+2072381165\,{t}^{19}+5859780674\,{t}^{18}-1127775119\,{t}^{17
}\\&&
-16287829166\,{t}^{16}-15833938922\,{t}^{15}+9251292427\,{t}^{14}\\&&
+38652814035\,{t}^{13}+44572754075\,{t}^{12}+10866248029\,{t}^{11}\\&&
-40129564125\,{t}^{10}-59007425756\,{t}^{9}-36829453004\,{t}^{8}\\&&
-10216139916\,{t}^{7}-63849664\,{t}^{6}+693364800\,{t}^{5}+187804368\,{
t}^{4}\\&&
\left. +24111840\,{t}^{3}+1694752\,{t}^{2}+63488\,t+1024 \right)
\end{array}$$

$$\begin{array}{lcl}c_1&=&
-768\,{t}^{31}+1488\,{t}^{30}+44636\,{t}^{29}-49538\,{t}^{28}-1198111
\,{t}^{27}\\&&
+359773\,{t}^{26}+19127286\,{t}^{25}+7602856\,{t}^{24}-
192894854\,{t}^{23}\\&&
-193322898\,{t}^{22}+1208988418\,{t}^{21}+
1968967542\,{t}^{20}\\&&
-4212884427\,{t}^{19}-10893520130\,{t}^{18}+
4099837581\,{t}^{17}\\&&
+31343794098\,{t}^{16}+22508418249\,{t}^{15}-
27437733598\,{t}^{14}\\&&
-67449042813\,{t}^{13}-62529644946\,{t}^{12}-
3629552721\,{t}^{11}\\&&
+84243589625\,{t}^{10}+130637976096\,{t}^{9}+
104165077688\,{t}^{8}\\&&
+50704667612\,{t}^{7}+15902199040\,{t}^{6}+
3290858704\,{t}^{5}\\&&
+451630576\,{t}^{4}+40498432\,{t}^{3}+2274080\,{t}^
{2}+72704\,t+1024
\end{array}$$

$$\begin{array}{lcl}c_2&=&
t\ \left( 680\,{t}^{29}-932\,{t}^{28}-39270\,{t}^{27}+33926\,{t}^{26}+1020385\,{t}^{25}\right.\\&&
-335552\,{t}^{24}-15580483\,{t}^{23}-2999586\,{t}^{22
}+151572589\,{t}^{21}\\&&
+104425399\,{t}^{20}-945694543\,{t}^{19}-
1131146667\,{t}^{18}\\&&
+3558106593\,{t}^{17}+6544185368\,{t}^{16}-
6226627151\,{t}^{15}\\&&
-20759382401\,{t}^{14}-4197042728\,{t}^{13}+29133893401\,{t}^{12}\\&&
+33005986439\,{t}^{11}+5921959164\,{t}^{10}-
22398414511\,{t}^{9}\\&&
-37477032816\,{t}^{8}-35648192872\,{t}^{7}-
21631096056\,{t}^{6}\\&&
-8298733544\,{t}^{5}-1992896768\,{t}^{4}-295849440
\,{t}^{3}\\&&
\left. -26179392\,{t}^{2}-1255424\,t-24832 \right)
\end{array}$$

$$\begin{array}{lcl}c_3&=&
- t^2\ \left( t-2 \right)  \left( t+1 \right)  \left( 2\,{t}^{4}+6\,{t}^{3}
-11\,{t}^{2}-30\,t-4 \right)\\&&  
\left( 124\,{t}^{21}-398\,{t}^{20}-5146
\,{t}^{19}+14694\,{t}^{18}+92616\,{t}^{17}-213234\,{t}^{16}\right. \\&&
-966327\,{t
}^{15}+1518831\,{t}^{14}+6391763\,{t}^{13}-5003278\,{t}^{12}\\&&
-26227554
\,{t}^{11}+1248286\,{t}^{10}+58532080\,{t}^{9}+36103178\,{t}^{8}\\&&
-41699603\,{t}^{7}-64544195\,{t}^{6}-41818519\,{t}^{5}-21472740\,{t}^{4
}\\&&\left. -8578026\,{t}^{3}-1961960\,{t}^{2}-218928\,t-9296 \right)
\end{array}$$
and 

$$\begin{array}{lcl}c_4&=&
({t}^{2}-2\,t-2 )  ( {t}^{2}-2\,t-7 ) ( t-2 ) ^{2} ( t+1 ) ^
{2}\\&& ( 2\,{t}^{4}+6\,{t}^{3}-11\,{t}^{2}-30\,t-4 ) ^{2}{t}^{
5}\\&&  ( 8\,{t}^{7}-10\,{t}^{6}-171\,{t}^{5}+209\,{t}^{4}+948\,{t}^{3}-721\,{t}
^{2}-1892\,t-249 ) 
\end{array}$$

The first coefficients of the series $y(t)$ are

$$-t+9\,{t}^{2}-49\,{t}^{3}+284
\,{t}^{4}-1735\,{t}^{5}+10955\,{t}^{6}-70695\,{t}^{7}+463087\,{t}^{8}$$
$$-3066450\,{t}^{9}+20471641\,{t}^{10}-137540539\,{t}^{11}+
928791019\,{t}^{12}\mp \dots$$

The \lq\lq complementary'' function 
$$\tilde y=\tilde g_\circ+\tilde g_N+\tilde g_{NW}+\tilde g_W+\tilde g_{SW}+
\tilde g_S+\tilde g_{SE}+\tilde g_E+\tilde g_{NE}$$
is associated to the matrice $J-M_1$ and $J-M_2$ (where $J$ is the square
matrix of order $9$ with all entries equal to $1$) satisfies of course
the same polynomial equation, after transposition of $y$ and $t$.

{\bf Remark.} For computations of huge terms in the series expansion of an
algebraic function, one
can use the following well-known trick: Any algebraic function
$y(t)=\sum a_n t^n$ of degre $d$ satisfies a linear differential equation 
$$\sum_{k=0}^d q_k(t) y^{(k)}(t)=0$$
with polynomial coefficients $q_0,\dots,q_t\in{\mathbf C}[t]$.
This allows a recursive computations of $a_n$
with time and memory requirements linear in $n$.
In our case, we get
$$q_0(t)y+q_1(t)y'+q_2(t)y''+q_3(t)y^{(3)}+q_4(t)y^{(4)}=0$$
where $q_0,\dots,q_4\in {\mathbf Z}[t]$ are polynomials of degrees
respectively $150,\ 151,\ 155,\ 156$ and $157$.

\subsection{Asymptotics}

The asymptotic growth rate of the coefficients of $y(t)$ 
is governed by the distance of the origin to the first ramification 
point of the corresponding sheet, see \cite {O}. Ramifications 
are above the zeros of the discriminant $D(t)$
of $P(y,t)$ with respect to $y$. This discriminant is given by

$$D(t)=r_1^2\ r_2\ r_3^2\ r_4^2\ r_\infty$$
where 
$$\begin{array}{lcl}
r_1&=&t^4-4\, t^3+6\, t^2+8\, t+1\\
r_2&=&
{t}^{13}-3\,{t}^{12}-16\,{t}^{11}+100\,{t}^{10}-86\,{t}^{9}-222\,{t}^{
8}+312\,{t}^{7}-544\,{t}^{6}\\
&&+4845\,{t}^{5}+10665\,{t}^{4}+9536\,{t}^{3
}+4084\,{t}^{2}+528\,t+16\\
r_3&=&
20\,{t}^{16}+50\,{t}^{15}-849\,{t}^{14}-937\,{t}^{13}+11563\,{t}^{12}+
7833\,{t}^{11}\\&&
-64177\,{t}^{10}-58882\,{t}^{9}+141152\,{t}^{8}+259280\,
{t}^{7}+82253\,{t}^{6}\\&&
-366913\,{t}^{5}-698955\,{t}^{4}-468324\,{t}^{3}
-122700\,{t}^{2}-13720\,t-552\\
r_4&=&5760\,{t}^{53}+8864\,{t}^{52}-425056\,{t}^{51}+\dots
-7200309248\,t-76152832\end{array}$$
($r_4$ is not involved in coarse asymptotics of $y(t)$).
The roots of the polynomial
$$r_\infty=t^6\ (t-2)^2\ (t+1)^3\ 
(2\,{t}^{4}+6\,{t}^{3}-11\,{t}^{2}-30\,t-4)^2$$
are critical points for the critical value $\infty$.

Since the coefficients of $y(t)$ have alternating signs, 
the \lq\lq smallest'' singularity of $y(t)$ is on the negative 
real halfline. The following table resumes
the relevant data for its computation. More precisely,
the algebraic function defined by $P$ has a ramification
of order $3$ with $y=\infty$ above $t=0$. The remaining sheet 
is unramified above $t=0$ and defines the generating 
function $y(t)$ under consideration.

The table contains the following informations: The first
column displays the argument $t$ considered in the corresponding
row. The second row indicates the factor of the discriminant $D(t)$
if $t$ is a root of $D(t)$. The remaining row displays information
about the inverse images of $t$ defined by the algebraic equation
of $y(t)$.

We have
$$\begin{array}{|l|c|l|}
\hline
t_0=0&r_\infty&y_1=y_2=y_3=\infty,\ y_4=0\\
\hline
t_0>t>t_1&&y_1<y_2<y_3<0<y_4\\
\hline
t_1\sim-0.04355&r_2&y_1\sim -2922<y_2=y_3\sim-1.083<0<y_4\sim0.066\\
\hline
t_1>t>t_2&&y_1<0<y_4,\ y_2=\overline{y_3}\in{\mathbf C}\setminus {\mathbf R}\\
\hline
t_2\sim-0.1118&r_3&y_1\sim-82.3<y_2=y_3\sim0.2194<y_4\sim0.4499\\
\hline
t_2>t>t_3&&y_1<0<y_4,\ y_2=\overline{y_3}\in{\mathbf C}\setminus {\mathbf R}\\
\hline
t_3\sim-0.14047& &0<y_4\sim 4.113,\ y_1=\infty,\ y_2=\overline{y_3}\sim
-2.98\pm 8.481i \\
\hline
t_3>t>t_4&&0<y_4<y_1,\ y_2=\overline{y_3}\in{\mathbf C}\setminus {\mathbf R}\\
\hline
t_4\sim-0.14118&r_\infty&0<y_4\sim8.692<y_1\sim28.07,\ y_2=y_3=\infty\\
\hline
t_4>t>t_5&&0<y_4<y_1,\ y_2=\overline{y_3}\in{\mathbf C}\setminus {\mathbf R}\\
\hline
t_5\sim-0.14127&r_2&0<y_4=y_1\sim14.89,\ y_2=\overline{y_3}\sim
23.95\pm 59.72i \\
\hline
\end{array}$$
The convergency radius of the series for $y(t)$ is of course given by
$\vert t_5\vert=-t_5$ and the asymptotic growth of the coefficients 
of $y(t)$ is roughly exponential with argument
$$1/t_5\sim -7.07857458512410303820641252737538586816317182$$

A slightly more precise asymptotical behaviour of the coefficients
of $y(t)$ can be computed as follows:

At the root 
$$\rho=t_5\sim -.14127137998962933757540882196178714222253950575630$$ 
of $r_2$, the ramified sheet corresponds
to the double root 
$$y_\rho\sim 14.88738808602894055277970788094544394$$ 
of $P(y,\rho)\in{\mathbf C}[y]$. (Caution: when computing
$y_\rho$ as a root of $P(y,\rho)$ one loses roughly 
half the digits since an error of order $\epsilon$ on $\rho$ induces an
error of order $\sqrt{\epsilon}$ on the corresponding two roots
approximating the double root $y_\rho$ of $P(y,\rho)$. A better
strategy is of course to compute $y_\rho$ as a (simple) root 
of the derivative $\frac{d}{dy}P(y,\rho)$ of $P(y,\rho)$).

In a open neighbourhood of $\rho$ we get now a Puiseux series expansion 
$$y(t)=h(t)+\sqrt{\rho-t}\ g(t)$$ with $h(t),\ g(t)$ holomorphic in an
open disc of radius 
$$.141478015962983913779403501>\rho\sim 
             .141271379989629337575408821$$ 
containing no other singularities of $y(t)$. 
The asymptotics of the generating function $y(t)$ 
are thus roughly given by
$$\begin{array}{cl}&\gamma_\rho\sqrt{\rho}\sqrt{1-t/\rho}\\
=&\gamma_\rho\sqrt{\rho}\sum_{n=0}^\infty {1/2\choose n} \left(
\frac{-t}{\rho}\right)^n\\
=&\gamma_\rho{\sqrt{\rho}}\left(1-\sum_{n=1}^\infty\frac{1}{2}\ 
\frac{1/2\cdot 3/2\cdots (n-3/2)}{n!}\left(\frac{t}{\rho}\right)^n\right)\\
=&\gamma_\rho{\sqrt{\rho}}\left(1-\sum_{n=1}^\infty
\frac{1\cdot 3\cdots (2n-3)}{2^n\ n!}\left(\frac{t}{\rho}\right)^n\right)\\
=&\gamma_\rho{\sqrt{\rho}}\left(1-\frac{1}{2}\sum_{n=1}^\infty
\frac{(2n-2)!}{4^{n-1}\  n!\ (n-1)!\ \rho^n}\ t^n\right)\end{array}$$
where $\gamma(\rho)=g(\rho)$. Since
$$\frac{(2n-2)!}{n!\ (n-1)!}\sim\frac{1}{n}\frac{\sqrt{4\pi(n-1)}
\ 4^{n-1}\ (n-1)^{2n-2}\ e^{-2n+2}}{2\pi\ (n-1)\ (n-1)^{2(n-1)}\ e^{-2(n-1)}}
\sim\frac{4^{n-1}}{\sqrt\pi\ n^{3/2}}$$
we get the asymptotics 
$$a_n\sim\frac{\gamma_\rho}{2\ \sqrt{\pi}\ n^{3/2}\ \rho^{n-1/2}}\ $$
The constant $\gamma_\rho$ can be computed by remarking that
$$\begin{array}{ccl}
0&=&P(h(t)+\sqrt{\rho-t}\ g(t),t)\\
&=&
P(y_\rho+\gamma_\rho\sqrt{\rho-t}+O((\rho-t)),t)\\
&=&\frac{\partial^2P}{\partial y^2}
\vert_{(y_\rho,\rho)}
\frac{\gamma_\rho^2\ (\rho-t)}{2}+\frac{\partial P}{\partial t}\vert_{
(y_\rho,\rho)}(t-\rho)+O((\rho-t)^{3/2}))
\end{array}$$
yielding 
$$\gamma_\rho\sqrt{\rho}=\sqrt{2\rho\frac{\frac{\partial P}{\partial t}\vert_{
(y_\rho,\rho)}}{\frac{\partial^2P}{\partial y^2}
\vert_{(y_\rho,\rho)}}}\sim \ 337.171657540870\ .$$
We have thus asymptotically
$$a_n\sim 95.11436852604511894068836\ \rho^{-n}\ n^{-3/2}$$
with $\rho\sim -.1412713799896293375754088219617871422225395057563006418$.
Unfortunately, the right hand side is a fairly correct approximation
of $a_n$ only for very huge values of $n$ (concretely, $n\sim 10^4$ yields
a only very few decimals). Indeed, the function $y(t)$ ramifies again for 
$$t\sim -.1414780159629839137794$$ 
which is extremely close to $\rho$.

Roland Bacher, INSTITUT FOURIER, Laboratoire de Math\'ematiques, 
UMR 5582 (UJF-CNRS), BP 74, 38402 St MARTIN  D'H\`ERES Cedex (France),
e-mail: Roland.Bacher@ujf-grenoble.fr

\end{document}